\documentclass[numbook, envcountsect, envcountsame, envcountreset, runningheads, smallextended]{svjour3}

\usepackage{graphicx}
\usepackage{fix-cm}
\usepackage{amssymb}
\usepackage{amsmath}
\usepackage{comment}
\usepackage[author-year,initials]{amsrefs}

\newtheorem{thm}{Theorem}[section]

\numberwithin{equation}{section}
\allowdisplaybreaks

\newtheorem{rem}{Remark}[section]
\newtheorem{defn}{Definition}[section]
\newtheorem{lem}{Lemma}[section]
\begin{document}

\title{Simulation and estimation for the fractional Yule process}


\author{Dexter O. Cahoy       \and
        Federico Polito
}


\institute{Dexter O. Cahoy $( \boxtimes )$ \at
       Program of Mathematics and Statistics\\
       College of Engineering and Science\\
		Louisiana Tech University, USA\\
        Tel: +1 318 257 3529\\
        Fax: +1 318 257 2182\\
         \email{dcahoy@latech.edu}
\and
Federico Polito \at
Dipartimento di Scienze Statistiche\\
 Sapienza University of Rome, Italy\\
 Tel: +39 0649910499\\
 Fax: +39 06 4959241\\
 \email{federico.polito@uniroma1.it}
}

\date{Received: 09-13-2010  / Accepted: date}

\maketitle

\begin{abstract}
        In this paper, we propose some representations of a generalized linear birth process called
		fractional Yule process (fYp). We also derive the probability distributions of the random birth and
		sojourn times. The inter-birth time distribution and the representations then yield algorithms on
		how to simulate	sample paths of the fYp. We also attempt to estimate the model parameters in order for the fYp to be usable in practice. The
		estimation procedure is then tested using simulated data as well. We also illustrate some major characteristics of  fYp which will be helpful for real applications.

\keywords{Yule--Furry process \and fractional calculus \and Mittag--Leffler \and Wright \and Poisson process \and birth process}
\subclass{	37A50 \and 62M86 \and 97K60}
\end{abstract}

\section{Introduction}

		The pure birth process is undoubtedly considered as one of the simplest branching processes.
		It has a Markovian structure and has already been extensively studied in the past. When the birth rate is linear, it is then usually called the pure linear birth or  classical Yule or Yule--Furry process (Yp). The pure linear birth process has been introduced by \ocite{mckendrick}, and  has been widely used to
		model various stochastic dynamical systems such as cosmic showers, epidemics, and population growth to name a few.  %In operations research,  a birth process was recently used by \ocite{dgt09}    %to study the valuation and risk analysis of multiname credit derivatives. In particular, a time-changed birth process was introduced to model correlated event timing in portfolios of %credit-sensitive securities such as  bonds and loans.

		For the sake of completeness, we review some known properties of the classical Yule process which will be used in the succeeding
		discussion. Let $\frak{N}(t)$ be the number of individuals in a Yule process with a single initial
		progenitor and birth intensity $\lambda > 0$. The $k$th state
		probability  or the probability of having exactly $k$ individuals  $\frak{p}_k(t) = \text{Pr} \left\{ \frak{N}(t) = k \mid \frak{N}(0)
		= 1 \right\}$ in a growing population at time $t>0$ solves the following Cauchy problem:
		\begin{equation}
			\label{diff}
			\begin{cases}
				\frac{d}{dt} \frak{p}_k(t) = - \lambda k \frak{p}_k(t) + \lambda (k-1)
				\frak{p}_{k-1} (t),
				& k \geq 1, \\
				\frak{p}_k(0) =
				\begin{cases}
					1, & k=1, \\
					0, & k > 1,
				\end{cases}
			\end{cases}
		\end{equation}
		where $\frak{p}_0(0) = 0$.  The explicit solution to \eqref{diff} is
		\[ %begin{equation}
			\frak{p}_k(t) = e^{- \lambda t} (1-e^{-\lambda t})^{k-1}, \qquad t>0, \: k \geq 1,
		\] %end{equation}
		with mean $\mathbf{E} \frak{N}(t)  = e^{\lambda t}$.
\begin{comment}
This model, however, turned out to be not satisfactory , as real phenomena display
		further complex properties, such as long memory and non-Markovian behaviour,
		which are neglected in the classical linear pure birth process.
\end{comment}
		To make the Yule process more flexible in taking into account more complex non-Markovian behaviour, some authors
		(\ocite{cahoy}, \ocite{pol}) proposed  a  more general model called the fractional Yule process (fYp).  A similar generalization of other point processes such as the Poisson process has previously been  carried out by \ocite{repin}, \ocite{jumarie}, \ocite{laskin}, \ocite{wang}, \ocite{mainardi},
		\ocite{wen}, \ocite{fan},  \ocite{mainardi2}, \ocite{cah}, \ocite{sibatov},
		\ocite{cahoy} and \ocite{orsbeg}.

		The aim of this paper is twofold: We want to derive related representations of fYp in terms of
		some classical or standard processes, and we want to construct algorithms on how to simulate a fYp
        and estimate the parameters.

		We organize the rest of the paper as follows: In  Section 2, we show  the fractional  generalization of the pure linear birth process.   In Section 3, it is illustrated that a pure linear birth process can also be viewed as a classical linear pure birth process with Wright-distributed 	random rates evaluated on a stretched time scale, i.e.,
		\[ %begin{equation}
			\mathfrak{N}^\nu(t) \overset{d}{=}\frak{N}_{\Xi}(t^\nu), \qquad \nu \in (0,1],
		\] %end{equation}
		where $\Xi$ is a random variable having the Wright probability density function
\begin{comment}
		\begin{equation}
			f_{\Xi}(\xi)=W_{-\nu, 1-\nu}(-\xi), \qquad \xi \in \mathbb{R}^+,
		\end{equation}
		with $W_{-\nu, 1-\nu}(-\xi)$ a specific case of the Wright function. It is defined as
\end{comment}
		\begin{equation} \label{bigxi}
			W_{-\nu, 1-\nu}(-\xi)=\sum_{r=0}^\infty \frac{(-\xi)^r}{r!\Gamma(1-\nu(r+1))}.
		\end{equation}		
Furthermore,  some Poisson-related representations are proved. In Section 4, we derive the birth and inter-birth time distributions. The structural representation,
fractional moments of the sojourn and birth times are also shown. In Section 5, we generate sample paths of a fYp  using our algorithms. In Section 6,  an estimation procedure is proposed using
		the moments of the log-transformed data, and  some empirical results are showed as well. Section 7 concludes
		the paper with a discussion on the key points and possible extensions of this study.

\begin{comment}
a  representation in terms of a mixed non-homogeneous
		Poisson process with intensity function $\lambda(t) = \Omega \lambda \Xi e^{\lambda \Xi t}$ is derived, where
		 $\Omega$ is a negative-exponential distributed random variable with mean 1.
		
    The following equalities in distribution:
		\begin{subequations}
			\begin{align} \notag
				& \mathfrak{N}^\nu(t) \overset{d}{=} M\left[ \Omega
				\left( e^{\lambda T_{2\nu}(t)}-1 \right) \right],
				\qquad t>0, \: \nu \in (0,1], \\ \notag
				& M (t) \overset{d}{=} \mathfrak{N}^\nu\left[\mathfrak{T}^\nu
				\left( \frac{1}{\lambda} \log
				\left( \frac{t}{\Omega}+1\right) \right) \right], \qquad t>0, \: \nu \in (0,1], \notag
			\end{align}
		\end{subequations}
		are also proved, where $M(t)$ is a homogeneous Poisson process with rate 1, and
		$\mathfrak{T}^\nu (t)$ is a process with the following one-dimensional distribution
		(given in terms of Fox functions):
		\[ %begin{equation}
			\text{Pr} \left\{ \mathfrak{T}^\nu(t) \in ds \right\}
			= t^{-\frac{1}{\nu}} H^{1,0}_{1,1} \left[ t^{-\frac{1}{\nu}}s
			\left|
			\begin{array}{l}
				\left(1-1/\nu,1/\nu\right) \\
				\left( 0,1 \right)
			\end{array}
			\right.
			\right] ds.
		\] %end{equation}
\end{comment}

\section{Generalization of the Yule process}
		
The fractional generalization of the Cauchy problem (\ref{diff}) was first carried out in
		\ocite{cahoy}, Section 8, and is described as follows:
		The authors defined the following difference-differential equations
		governing the state probabilities $p_k^\nu (t) = \text{Pr} \left\{ N^\nu (t) = k \mid
		N^\nu (0) =1 \right\}$:
		\begin{equation}
			\label{gov}
			\frac{\mathfrak{d}^\nu}{\mathfrak{d}t^\nu} p_k^\nu(t) = \lambda \left[
			\sum_{l=1}^{k-1} p_l^\nu(t) p_{k-l}^\nu (t) - p_k^\nu(t) \right] + \frac{t^{-\nu}}{\Gamma (1-\nu)}
			\delta_{k,1}, \qquad \nu \in (0,1], \: k \geq 1,
		\end{equation}
		where the initial condition
		\[%begin{equation}
			p_k^\nu(0) =
			\begin{cases}
				1, & k=1, \\
				0, & k>1,
			\end{cases}
		\]%end{equation}
		is incorporated into equation \eqref{gov} through the Kronecker delta $\delta_{k,1}$.
		The fractional derivative appearing in \eqref{gov} is the so-called Riemann--Liouville
		operator, and	is defined as
		\begin{equation}
			\label{rl}
			\begin{cases}
				\frac{\mathfrak{d}^\nu}{\mathfrak{d} t^\nu} f (t) =
				\frac{1}{\Gamma (1-\nu)} \frac{d}{dt} \int_0^t \frac{f(s)}{(t-s)^\nu} ds,
				& \nu \in (0,1), \\
				f' (t), & \nu = 1.
			\end{cases}
		\end{equation}
		Furthermore, the mean number of individuals in the system was found to be
		\begin{equation} \label{mun}
			\mathbf{E} \bigl[N^\nu (t) \bigr] =
			E_{\nu,1} (\lambda t^\nu), \qquad t>0, \: \nu \in (0,1],
		\end{equation}
		where
		\[ %begin{equation}
			E_{\alpha, \beta} (\tau) = \sum_{r=0}^\infty \frac{\tau^r}{\Gamma (\alpha r+\beta)},
			\qquad \alpha, \beta, \in \mathbb{R}^+,\: \tau \in \mathbb{R},
		\] %end{equation}
		is the generalized Mittag--Leffler function.
\begin{comment}
Its Laplace transform is
		\begin{equation}
			\int_0^\infty e^{-zt} E_{\alpha,\beta} (- \vartheta t^\alpha) dt =
			\frac{z^{\alpha -\beta}}{
			z^\alpha + \vartheta}.
		\end{equation}
\end{comment}

	    Let  $\mathfrak{N}^\nu
		(t)$ be the number of individuals in a fractional linear
		birth process or fractional Yule or Yule--Furry process (fYp) up to the time $t >0$. The
		state probabilities $\mathfrak{p}_k^\nu (t) = \text{Pr} \left\{ \mathfrak{N}^\nu (t) = k
		\mid \mathfrak{N}^\nu (0) = 1 \right\}$ solve the
		following Cauchy problem:
		\begin{equation}
			\label{diff-frac}
			\begin{cases}
				\frac{d^\nu}{dt^\nu} \mathfrak{p}_k^\nu(t) =
				- \lambda k \mathfrak{p}_k^\nu(t) + \lambda (k-1) \mathfrak{p}_{k-1}^\nu (t),
				&  k \geq 1, \\
				\mathfrak{p}_k^\nu(0) =
				\begin{cases}
					1, & k=1, \\
					0, & k > 1,
				\end{cases}
			\end{cases}
		\end{equation}
		which is also a fractional generalization of \eqref{diff}. The fractional derivative involved in
		\eqref{diff-frac} is  now the Caputo operator,  and is defined as
		\begin{equation}
			\label{caputo}
			\begin{cases}
				\frac{d^\nu}{d t^\nu} f (t) =
				\frac{1}{\Gamma (1-\nu)} \int_0^t \frac{f'(s)}{(t-s)^\nu} ds,
				& \nu \in (0,1), \\
				f' (t), & \nu = 1.
			\end{cases}
		\end{equation}

		Moreover, the Riemann--Liouville \eqref{rl} and the Caputo \eqref{caputo} fractional
		derivatives are linked together by the following relation (see
		\ocite{kilbas}, page 91):
		\begin{equation}
			\label{relation}
			\frac{d^\nu}{d t^\nu} f (t) = \frac{\mathfrak{d}^\nu}{\mathfrak{d} t^\nu} f (t)
			-\frac{f(0)}{\Gamma (1-\nu)} t^{-\nu}, \qquad \nu \in (0,1).
		\end{equation}
        From \eqref{relation}, it is easy to see that both fractional derivatives coincide
		when $f(0)=0$ for each $k>1$.  The solution to the Cauchy problem
		\eqref{diff-frac} is
		\begin{equation}
			\label{law}
			\mathfrak{p}_k^\nu (t)=\sum_{l=1}^k \binom{k-1}{l-1}(-1)^{l-1}
			E_{\nu,1}(-\lambda l t^\nu), \qquad
			k\geq1, \: \nu \in (0,1].
		\end{equation}

Note that the mean number of individuals $\mathbf{E} \bigl[ \mathfrak{N}^\nu (t) \bigr] $ in the fractional Yule process is the same as (\ref{mun}), and the variance can be calculated as
\[
\mathbf{V}ar \left( \mathfrak{N}^\nu (t) \right) = 2 E_{\nu,1}\left(2\lambda  t^\nu \right)-E_{\nu,1}\left(\lambda  t^\nu \right)-E_{\nu,1}^2\left(\lambda t^\nu \right)
\]
From here on, we emphasize that the fractional derivative operation is performed in Caputo's sense. 	
\begin{comment}
We now state a useful subordination relation for the fractional linear birth process
		which permits us to derive further representations of fYp. The fractional Yule process
		$\mathfrak{N}^\nu(t)$ can be considered as a classical Yule process evaluated at
		a random time, i.e.\
		\begin{equation}
			\label{sub}
			\mathfrak{N}^\nu(t) \overset{d}{=} \frak{N}(T_{2\nu}(t)), \qquad t>0, \: \nu \in (0,1],
		\end{equation}
		where $T_{2\nu}(t)$ has one-dimensional distribution that solves the following fractional diffusion equation:
		\begin{equation}
			\label{diffusion}
			\begin{cases}
				\frac{\partial^{2\nu}}{\partial t^{2\nu}} g (x,t) =
				\frac{\partial^2}{\partial x^2} g(x,t),
				& t>0,  \:\nu \in (0,1],  \: x>0,\\
				\left. \frac{\partial}{\partial x}g(x,t) \right|_{x=0}=0,\\
				g(x,0)=\delta(x),
			\end{cases}
		\end{equation}
		with the initial condition $g_t(x,0)=0$, when $1/2<\nu \leq 1$. Details about subordination using $T_{2 \nu} (t)$
are provided in  \ocite{ors2008} and  \ocite{pol}.
\end{comment}

	\section{Stretched Yule process with random rates and related representations}

		In this section, we present some relevant and interesting representations  of the fractional Yule process (fYp).
		We start by proving a subordination relation that links the fractional Yule process with its classical
		counterpart. %We report the Theorem here for the sake of %completeness.
		\begin{thm}
			\label{thea}
			Let $\mathfrak{N}^\nu(t)$ be the number of individuals in  a fractional
			Yule process at time $t>0$. Then the following equality in distribution holds:
			\begin{equation}
				\label{sub2}
				\mathfrak{N}^\nu(t) \overset{d}{=}  \frak{N}(T_{2\nu}(t)),
			\end{equation}
			where $\frak{N}(t)$ is a classical Yule process,  $\nu \in (0,1]$, and
			$T_{2\nu}(t)$  is a random time whose distribution coincides with
			the solution to the following Cauchy problem
			\begin{equation}
				\label{diffusion2}
				\begin{cases}
					\frac{\partial^{2\nu}}{\partial t^{2\nu}} g (x,t) =
					\frac{\partial^2}{\partial x^2} g(x,t),
					& x>0,\\
					\left. \frac{\partial}{\partial x}g(x,t) \right|_{x=0}=0,\\
					g(x,0)=\delta(x),
				\end{cases}
			\end{equation}
			with the initial condition $g_t(x,0)=0$, when $1/2<\nu \leq 1$.

			\begin{proof}
				Let $G^\nu(u,t)$, $t>0$, $|u|<1$, be the probability generating function
				of the fractional Yule process. To prove \eqref{sub2}, 	it is sufficient to observe that
				\begin{align} \notag
					& \int_0^\infty e^{-zt} G^\nu(u,t) dt \\
					& = \int_0^\infty e^{-zt} \sum_{k=1}^\infty
					u^k \mathfrak{p}_k^\nu (t) dt \notag \\
					& = \int_0^\infty e^{-zt} \sum_{k=1}^\infty u^k
					\sum_{l=1}^k \binom{k-1}{l-1} (-1)^{l-1} E_{\nu,1}(-\lambda l t^\nu)
					dt \notag \\
					& = \sum_{k=1}^\infty u^k \sum_{l=1}^k \binom{k-1}{l-1}
					(-1)^{l-1} \frac{z^{\nu-1}}{z^\nu+\lambda l} \notag \\
					& = \sum_{k=1}^\infty u^k \sum_{l=1}^k \binom{k-1}{l-1}
					(-1)^{l-1} z^{\nu-1} \int_0^\infty e^{-s(\lambda l+z^\nu)} ds \notag \\
					& = \int_0^\infty \sum_{k=1}^\infty u^k \sum_{l=1}^k
					\binom{k-1}{l-1} (-1)^{l-1} e^{-s \lambda l} z^{\nu-1}
					e^{-sz^\nu} ds \notag \\
					& = \int_0 ^\infty \sum_{k=1}^\infty u^k \sum_{l=1}^k \binom{k-1}{l-1}
					(-1)^{l-1} e^{-\lambda l s} \int_0^\infty e^{-zt}
					\text{Pr} \left\{ T_{2\nu}(t) \in ds \right\} dt \notag \\
					& = \int_0^\infty e^{-zt} \left[
					\sum_{k=1}^\infty u^k \int_0^\infty \text{Pr} \left\{ \mathfrak{N}(s) = k \right\}
					\text{Pr} \left\{ T_{2 \nu} (t) \in ds \right\} \right] dt \notag \\
					& = \int_0^\infty e^{-zt} \left[ \sum_{k=1}^\infty u^k
					\text{Pr} \left\{ \mathfrak{N}(T_{2 \nu}(t) ) = k \right\} \right] dt. \quad \qed \notag
				\end{align}
			\end{proof}
		\end{thm}
		
		\begin{rem}
			Note that, the solution to \eqref{diffusion2}, also solves the fractional differential equation
			\begin{align}
				\frac{\partial^\nu}{\partial t^\nu} g(x,t) = - \frac{\partial}{\partial x} g(x,t).
			\end{align}
		\end{rem}		
		
		\begin{rem}
			In the proof of Theorem \ref{thea}, we used the Laplace transform
			of $\text{Pr} \left\{ T_{2\nu}(t) \in ds \right\}$ which is
			\[ %begin{equation}
				\int_0^\infty e^{-zt} \text{Pr} \left\{ T_{2\nu}(t) \in ds \right\} =
				z^{\nu-1} e^{-sz^\nu} ds, \qquad s>0.
			\] %end{equation}
		\end{rem}

		In the next Theorem, we derive a random-rate representation of the fractional Yule
		process using the preceding subordination relation.
		\begin{thm}{(Representation A)}
			\label{repra} Let $t>0$ and $\nu \in (0,1]$. Then the following equality in distribution holds:
			\begin{equation}
				\label{random}
				\mathfrak{N}^\nu(t) \overset{d}{=} \frak{N}_\Xi(t^\nu),
			\end{equation}
			where $\frak{N}_\Xi(t^\nu)$ is a classical linear birth process with random
			rate $\lambda \Xi$ evaluated at $t^\nu$, $\Xi$ is a Wright-distributed
			random variable with probability density function $W_{-\nu, 1-\nu}(-\xi)$ in (\ref{bigxi}).
\begin{comment}
			\begin{equation}
				\label{bigxi}
				f_{\Xi}(\xi)=W_{-\nu, 1-\nu}(-\xi), \qquad \xi \in \mathbb{R}^+,
				\end{equation}
				and
				\begin{equation}
				W_{-\nu, 1-\nu}(-\xi)=\sum_{r=0}^\infty \frac{(-\xi)^r}{r!(1-\nu(r+1))}.
			\end{equation}
\end{comment}
			\begin{proof}
				To prove equality \eqref{random}, we use the
				subordination relation \eqref{sub2} as follows:
				\begin{align}
					\label{formula}
					& \text{Pr} \left\{ \mathfrak{N}^\nu(t) = k \mid \mathfrak{N}^\nu(0)=1
					\right\} \\
					& = \int_0^\infty \text{Pr} \left\{ \frak{N}(s) = k \mid
					\frak{N}(0)=1 \right\}
					\text{Pr} \left\{ T_{2\nu}(t) \in ds \right\} \notag \\
					& = \int_0^\infty \sum_{l=1}^k \binom{k-1}{l-1} (-1)^{l-1}
					e^{-\lambda ls} t^{-\nu} W_{-\nu,1-\nu} (-t^{-\nu} s) ds \notag \\
					& = \int_0^\infty \sum_{l=1}^k \binom{k-1}{l-1} (-1)^{l-1}
					e^{-\lambda l \xi t^\nu} W_{-\nu,1-\nu} (-\xi) d \xi \notag \\
					& = \int_0^\infty \text{Pr} \left\{ \frak{N}_\xi (t^\nu) = k
					\mid \frak{N}_\xi (0) = 1 \right\} W_{-\nu,1-\nu}(-\xi) d \xi, \notag
				\end{align}
				and this leads to \eqref{random}. \qed
			\end{proof}
		\end{thm}
		Note that in the second step of formula \eqref{formula}, we used the explicit form
		of the solution to the fractional diffusion equation \eqref{diffusion2} which is
		(see \ocite{podlubny}, formula (4.22), page 142)
		\[ %begin{equation}
			\text{Pr} \left\{ T_{2 \nu} (t) \in ds \right\} = t^{-\nu}
			W_{-\nu,1-\nu} (-t^{-\nu} s) ds, \qquad s>0.
		\] %end{equation}
		\begin{rem}
			As noted above, representation \eqref{random} holds for the one-dimensional
			state probability distribution $\mathfrak{p}_k^\nu (t)$, $t>0$,
			$k \geq 1$. This, however is sufficient in the sense that the process
			$\frak{N}_\Xi(t^\nu)$  has distribution that solves \eqref{diff-frac}.
		\end{rem}

		We now prove a further interesting representation of the fractional Yule process in
		terms of a specific mixed non-homogeneous Poisson process.

		Starting from the second-to-last step of formula \eqref{formula}, we obtain
		\begin{align} \notag
			\mathfrak{p}^\nu_k (t) & = \int_0^\infty \sum_{l=1}^k \binom{k-1}{l-1}
			(-1)^{l-1} e^{-\lambda l \xi t^\nu} W_{-\nu,1-\nu}(-\xi) d\xi \\ \notag
			& = \int_0^\infty e^{-\lambda \xi t^\nu} \left[1-e^{-\lambda \xi t^\nu}\right]^{k-1}
			W_{-\nu,1-\nu}(-\xi) d\xi \notag \\
			& = \int_0^\infty \frac{1}{\left[e^{\lambda \xi t^\nu}\right]^{k}}
			\left[ e^{\lambda \xi t^\nu} -1 \right]^{k-1} W_{-\nu,1-\nu}(-\xi) d\xi. \notag
		\end{align}
		Recalling the identity
		\[ %begin{equation}
			\int_0^\infty e^{-ax} x^r dx = a^{-(r+1)} r!, \qquad r \in \mathbb{N},\:
			\mathfrak{R}(a)>0,
		\]%end{equation}
		we get
		\begin{align}
			\label{silent}
			\mathfrak{p}_k^\nu(t) & = \int_0^\infty \int_0^\infty e^{-\omega e^{\lambda
			\xi t^\nu}} \omega^{k-1} \frac{\left[ e^{\lambda \xi t^\nu}-1 \right]^{k-1}}{(k-1)!}
			W_{-\nu,1-\nu}(-\xi) d\omega \, d\xi \\
			& = \int_0^\infty \int_0^\infty \frac{e^{-\omega \left[ e^{\lambda \xi t^\nu}-1
			\right]} \omega^{k-1} \left[ e^{\lambda \xi t^\nu}-1 \right]^{k-1} }{(k-1)!}
			e^{-\omega} W_{-\nu,1-\nu}(-\xi) d \omega \, d\xi \notag \\
			& = \int_0^\infty \int_0^\infty \frac{ e^{-\int_0^{t^\nu} \omega \lambda \xi
			e^{\lambda \xi s} ds} \left[ \int_0^{t^\nu} \omega \lambda \xi
			e^{\lambda \xi s} ds \right]^{k-1} }{(k-1)!}
			e^{-\omega} W_{-\nu,1-\nu}(-\xi) d\omega \, d\xi.
			\notag
		\end{align}
		Thus, we have obtained a representation in terms of a mixed non-homogeneous
		Poisson process with intensity function
		\[ %begin{equation}
			\lambda(t)=\Omega \lambda \Xi e^{\lambda \Xi t}, \qquad t>0,
		\] %end{equation}
		where the distribution of $\Omega$ is negative-exponential with mean equal to 1, and
		$\Xi$ has probability density function \eqref{bigxi}. Note that the random variable
		$\Omega$, conditional on $\Xi=\xi$, is such that
		\[ %begin{equation}
			\frac{N_\xi (t^\nu)}{\mathbf{E} N_\xi (t^\nu)} \overset{\text{a.s.}}{\longrightarrow}
			\Omega,
		\] %end{equation}
		as $t\rightarrow \infty$ (see e.g.\ \ocite{keiding}, \ocite{waugh}, \ocite{harris}).

		\begin{rem}
			A simple change of variable also allows us to obtain a representation
			in terms of a mixed non-homogeneous Poisson process evaluated at the random
			time $T_{2\nu}(t)$, $t>0$. From the second step of formula
			\eqref{silent}, we have
			\begin{align} \notag
				\mathfrak{p}_k^\nu(t) & =
				\int_0^\infty \int_0^\infty \frac{e^{-\omega \left[ e^{\lambda \xi t^\nu}-1
				\right]} \omega^{k-1} \left[ e^{\lambda \xi t^\nu}-1 \right]^{k-1} }{(k-1)!}
				e^{-\omega} W_{-\nu,1-\nu}(-\xi) d \omega \, d\xi \\ \notag
				& = \int_0^\infty \int_0^\infty
				\frac{e^{-\omega \left[ e^{\lambda s}-1
				\right]} \omega^{k-1} \left[ e^{\lambda s}-1 \right]^{k-1} }{(k-1)!}
				e^{-\omega} \frac{1}{t^\nu} W_{-\nu,1-\nu} \left( -\frac{s}{t^\nu} \right)
				ds \, d\omega. \notag
			\end{align}
			Consider a non-homogeneous Poisson process $\mathrm{N}(t)$
			with intensity function $\lambda(t) = \Omega \lambda e^{\lambda t}$. Then the state probabilities of the fractional Yule process can be written as
			\begin{align}
				\label{conditional}
				\mathfrak{p}_k^\nu(t) & = \int_0^\infty e^{-\omega} \int_0^\infty
				\text{Pr} \left\{ \mathrm{N}(s) = k-1 \right\}
				\text{Pr} \left\{ T_{2 \nu}(t) \in ds \right\} d \omega \\
				& = \mathbf{E}_\Omega \mathrm{N} (T_{2\nu}(t)). \notag
			\end{align}
		\end{rem}

		In addition, the subordinated non-homogeneous Poisson process $\mathrm{N}(T_{2\nu}(t))$ conditioned on $\Omega = \omega$ could be interesting as the fractional homogeneous Poisson process admits a similar representation \cite{orsbeg}.

		Let $q_k^\nu(t)$  be the state probabilities of
		$\mathrm{N} (T_{2\nu} (t))$, i.e.,\
		\[ %begin{equation}
			q_k^\nu(t) = \text{Pr} \left\{ \mathrm{N} (T_{2\nu}(t)) = k-1
			\right\}, \qquad t>0, \: k \geq 1.
		\] %end{equation}
		Then
		\begin{equation}
			\label{study}
			q_k^\nu(t) = \int_0^\infty \frac{ e^{-\omega \left[ e^{\lambda s}-1 \right]}
			\omega^{k-1} \left[ e^{\lambda s}-1 \right]^{k-1} }{(k-1)!}
			\text{Pr} \left\{ T_{2\nu}(t) \in
			ds \right\}.
		\end{equation}
		Applying the Laplace transform to \eqref{study}, we have
		\begin{align} \notag
			\int_0^\infty e^{-zt} q_k^\nu(t) dt & =
			\int_0^\infty \frac{ e^{-\omega \left[ e^{\lambda s}-1 \right]}
			\omega^{k-1} \left[ e^{\lambda s}-1 \right]^{k-1} }{(k-1)!} z^{\nu-1} e^{-sz^\nu} ds \\ \notag
			& = \int_0^\infty e^\omega \frac{ e^{-\omega e^{\lambda s}} \omega^{k-1}
			\left[ 1-e^{\lambda s} \right]^{k-1} }{(k-1)!} (-1)^{k-1} z^{\nu-1} e^{-sz^\nu} ds, \notag
		\end{align}
		and by taking into account the relations
		\[ %begin{equation}
			e^{-\omega e^{\lambda s}} = \sum_{l=0}^\infty \frac{(-\omega)^l e^{\lambda sl}}{l!},
		\] %end{equation}
		\[ %begin{equation}
			\left[ 1-e^{\lambda s} \right]^{k-1} = \sum_{j=0}^{k-1} \binom{k-1}{j} (-1)^j
			e^{\lambda sj},
		\] %end{equation}
		we arrive at the equality
		\begin{align}
			\label{floor}
			& \int_0^\infty e^{-zt} q_k^\nu(t) dt \\
			& = \int_0^\infty \frac{e^\omega}{(k-1)!} (-1)^{k-1} \omega^{k-1} \sum_{l=0}^\infty
			\sum_{j=0}^{k-1} \frac{(-\omega)^l}{l!} e^{\lambda ls} \binom{k-1}{j}
			(-1)^j e^{\lambda sj} z^{\nu-1} e^{-sz^\nu} ds \notag \\
			& = \frac{e^\omega}{(k-1)!} (-1)^{k-1} \omega^{k-1}
			\sum_{l=0}^\infty \frac{(-\omega)^l}{l!}
			\sum_{j=0}^{k-1} \binom{k-1}{j} (-1)^j z^{\nu-1} \int_0^\infty e^{-s \left[
			z^\nu -\lambda \left( l+j \right) \right]} ds \notag \\
			& = \frac{e^\omega}{(k-1)!} (-\omega)^{k-1} \sum_{l=0}^\infty \frac{(-\omega)^l}{l!}
			\sum_{j=0}^{k-1} \binom{k-1}{j} (-1)^j \frac{z^{\nu-1}}{z^\nu-\lambda \left(l+j\right)}.
			\notag
		\end{align}
		Applying the inverse Laplace transform to equation \eqref{floor}, we obtain the explicit
		expression of the state probabilities as
		\begin{equation}
			\label{url}
			q_k^\nu(t) = \frac{e^\omega (-\omega)^{k-1}}{(k-1)!} \sum_{l=0}^\infty
			\frac{(-\omega)^l}{l!}
			\sum_{j=1}^k \binom{k-1}{j-1} (-1)^{j-1} E_{\nu,1} \left[
			\lambda \left( l+j \right) t^\nu \right], \quad k \geq 1.
		\end{equation}

		\begin{rem} From equation \eqref{url}, it is straightforward to obtain  the classical form of the
			state probabilities of the (conditional) non-homogeneous
			Poisson process ($\nu=1$) with intensity function $\lambda(t) = \omega \lambda
			e^{\lambda t}, t>0$.
		\end{rem}

		We introduce a definition and a	lemma below which will be helpful in transforming fYp into a
		non-homogeneous Poisson process with  rate 1. In order to do so, we present here the standard
		definition, by means of a Mellin--Barnes type integral, of the so-called Fox function:
		\begin{align}
			H_{p,q}^{m,n} \left[ x \left|
			\begin{array}{l}
				(a_1,A_1),\dots,(a_p,A_p) \\
				(b_1,B_1),\dots,(b_p,B_p)
			\end{array}
			\right. \right] = \frac{1}{2\pi i} \int_{c -i\infty}^{c +i\infty}
			\Theta(z)x^{-z} dz, \quad x \neq 0,
		\end{align}
		where
		\begin{align}
			\Theta (z) = \frac{\left\{ \prod_{j=1}^m \Gamma(b_j+B_j z) \right\}
			\left\{ \prod_{j=1}^n \Gamma(1-a_j - A_j z) \right\}}{\left\{ \prod_{j=m+1}^q
			\Gamma (1-b_j-B_j z) \right\}\left\{ \prod_{j=n+1}^p \Gamma(a_j+A_j z) \right\}}.
		\end{align}
		Each empty product is interpreted as unity. For more information on Fox functions
		we refer to \ocite{mathai}.

		\begin{defn}
			Let $\mathfrak{T}^\nu(t)$ be a random time process
			whose one-dimensional distribution is given by
			\[ %begin{equation}
				\text{Pr} \left\{ \mathfrak{T}^\nu(t) \in ds \right\}
				= h(t,s) ds = t^{-\frac{1}{\nu}} H^{1,0}_{1,1} \left[ t^{-\frac{1}{\nu}}s
				\left|
				\begin{array}{l}
					\left(1-1/\nu,1/\nu\right) \\
					\left( 0,1 \right)
				\end{array}
				\right.
				\right] ds,
			\] %end{equation}
			where $t>0$, $s>0$,  $\nu \in (0,1]$.
			Furthermore, $h(t,s)$ has Mellin transform
			\begin{equation}
				\label{mellin}
				\int_0^\infty s^{\eta-1} h(t,s) ds = \frac{\Gamma (\eta)}{
				\Gamma\left(1-\frac{1}{\nu}+\frac{1}{\nu}\eta\right)} t^{\frac{\eta-1}{\nu}}.
			\end{equation}
		\end{defn}

		\begin{lem}
			\label{lemma-inv}
			Let $\mathfrak{N}^\nu(t)$ be a fractional Yule process with rate
			$\lambda>0$ and $t>0$. Then the process $\mathfrak{N}^\nu(\mathfrak{T}^\nu(t))$  is a classical
			Yule process with rate $\lambda$.
			\begin{proof}
				Define $G^\nu(u,t)$ and $G(u,t)$, $t>0$, $|u|\leq 1$ as the probability
				generating functions of fYp and
				the classical Yule process,  respectively. Then
				\[ %begin{equation}
					\int_0^\infty G^\nu(u,s) h(t,s) ds =
					\int_0^\infty \sum_{k=1}^\infty u^k \sum_{j=1}^k
					\binom{k-1}{j-1} (-1)^{j-1} E_{\nu,1} (-\lambda j s^\nu)
					h(t,s) ds.
				\] %end{equation}
				In the following we use the Mellin--Barnes representation of the Mittag--Leffler
				function
				\[ %begin{equation}
					E_{\nu,1}(x) = \frac{1}{2\pi i}
					\int_{c -i\infty}^{c +i\infty}
					\frac{\Gamma(z) \Gamma(1-z)}{\Gamma(1-\nu z)} (-x)^{-z} dz, \qquad \nu > 0,\: x \neq 0
				\] %end{equation}
				(see \ocite{kilbas}, page 41, formula (1.8.14)). Note that when $\nu=1$ we retrieve
				the Mellin--Barnes representation of the exponential function
				\begin{align}
					e^x = \frac{1}{2\pi i}
					\int_{c -i\infty}^{c +i\infty}
					\Gamma(z) (-x)^{-z} dz, \qquad x \neq 0.
				\end{align}
				(see \ocite{paris}, page 89, formula (3.3.2)).
				
				We obtain
				\begin{align} \notag
					& \int_0^\infty G^\nu (u,s) h(t,s) ds \\ \notag
					& = \sum_{k=1}^\infty u^k \sum_{j=1}^k \binom{k-1}{j-1}
					\frac{(-1)^{j-1}}{2 \pi i} \int_{c -i \infty}^{c+i \infty}
					\frac{\Gamma(z) \Gamma(1-z)}{\Gamma(1-\nu z)} (\lambda j)^{-z}
					\int_0^\infty \frac{h(t,s)}{s^{\nu z}} ds\, dz. \notag
				\end{align}
				Applying formula \eqref{mellin}, we can write
				\begin{align}
					\label{oxford}
					& \int_0^\infty G^\nu (u,s) h(t,s) ds \\
					& = \sum_{k=1}^\infty u^k \sum_{j=1}^k \binom{k-1}{j-1}
					\frac{(-1)^{j-1}}{2 \pi i} \int_{c -i \infty}^{c+i \infty}
					\frac{\Gamma(z) \Gamma(1-z)}{\Gamma(1-\nu z)} (\lambda j)^{-z}
					\frac{\Gamma (1-\nu z)}{\Gamma (1-z)} t^{-z} dz \notag \\
					& = \sum_{k=1}^\infty u^k \sum_{j=1}^k \binom{k-1}{j-1}
					\frac{(-1)^{j-1}}{2 \pi i} \int_{c -i \infty}^{c+i \infty}
					\Gamma(z) (\lambda j t)^{-z} dz \notag \\
					& = \sum_{k=1}^\infty u^k \sum_{j=1}^k \binom{k-1}{j-1}
					(-1)^{j-1} e^{-\lambda j t} \notag \\
					& = \sum_{k=1}^\infty u^k e^{-\lambda t} \left[ 1-e^{-\lambda t}
					\right]^{k-1} \notag \\
					& = G(u,t). \qquad \qed \notag
				\end{align}
			\end{proof}
		\end{lem}

		\begin{rem}
			Note that it is straightforward to generalize Lemma \ref{lemma-inv} to the
			more general (non-linear) case.
		\end{rem}

		\begin{rem}
			Letting $u=1$ in \eqref{oxford}, we have
			\begin{align} \notag
				& \sum_{k=1}^\infty \int_0^\infty \mathfrak{p}^\nu_k(s) h(t,s) ds =
				\sum_{k=1}^\infty \frak{p}_k(t) \\ \notag
				& \Leftrightarrow \int_0^\infty h(t,s) ds = 1. \notag
			\end{align}
		\end{rem}

		\begin{thm}
			Consider a fractional Yule process $\mathfrak{N}^\nu(t)$
			with birth rate $\lambda>0$, $t>0$, and $\nu \in (0,1]$. 	Then the random time-changed process
			\[ %begin{equation}
				\mathfrak{N}^\nu\left[\mathfrak{T}^\nu \left( \frac{1}{\lambda} \log
				\left( \frac{t}{\Omega}+1\right) \right) \right]
			\] %end{equation}
			has one-dimensional distribution which coincides with that of a non-homogeneous
			Poisson process $M(t)$ with rate 1.
			\begin{proof}
				It readily follows from \eqref{conditional}, Lemma \ref{lemma-inv}
				and Theorem 1 of \ocite{kendall}.
			\end{proof}
		\end{thm}

	\section{Wait and sojourn time distributions}

		We now show that the sojourn or inter-birth time of fYp follows the Mittag--Leffler distribution.
		Let $T_i^\nu$, $i \geq 1$, denote the time between the  $(i-1)$th and $i$th birth.
		This means that $T_i^\nu$ is the time it takes for the population size to grow from $i$ to $i+1$.
		More specifically, we will show that the sojourn times $T_i^\nu$'s are independent and
		$T_i^\nu$ is distributed
		\begin{equation}
			\label{interbirth}
			f_{T_i^\nu} (t) = i \lambda t^{\nu-1} E_{\nu, \nu} (-i \lambda t^\nu ), \quad i \geq 1.
		\end{equation}
		Recall that when $\nu=1$, the inter-birth times $T_i$'s of the Yp are independent and $T_i$ is
		exponentially distributed with rate $i \lambda$, $i \geq 1$.  Moreover, the waiting or birth time distribution for the pure linear birth
        process ($\nu=1$) satisfies the following two equalities:
		%\begin{subequations}
			\begin{align} \notag
				&\text{Pr} ( \mathfrak{W}_j=  T_1 + \cdots + T_j \leq t) =
				\text{Pr} (\mathfrak{N}(t) \geq j+1 \big|  \mathfrak{N}(0) = 1 ) \\ \notag
				\text{and} & \\ \notag
				&\frak{p}_j (t) = \text{Pr} ( \mathfrak{W}_{j-1} \leq t)
				- \text{Pr} ( \mathfrak{W}_j \leq t).
			\end{align}
		%\end{subequations}

Let $\mathfrak{W}_j^\nu = T_1^\nu + T_2^\nu + \cdots + T_j^\nu$ be the waiting
		time of the $j$th birth of the fYp.  We now show that the preceding two equations hold true as well for the
        fractional or general case $(0< \nu \leq 1)$, i.e.,
		\begin{subequations}
			\begin{align}
				& \label{eq2a} \text{Pr} ( \mathfrak{W}_j^\nu \leq t) =
				\text{Pr} (\mathfrak{N}^\nu(t) \geq j+1 \big|
				\mathfrak{N}^\nu(0) = 1 ), \quad j\geq 1, \\
				\text{and} & \notag \\
				& \label{eq2b} \frak{p}_j^\nu(t) = \text{Pr} ( \mathfrak{W}_{j-1}^\nu \leq t)
				- \text{Pr} ( \mathfrak{W}_j^\nu \leq t).
			\end{align}
		\end{subequations}
        Using \eqref{law}, we obtain
		\begin{align}
			\label{eqprob}
			\text{Pr} (\mathfrak{N}^\nu(t) \geq j+1 \big|  \mathfrak{N}^\nu(0) = 1 )  &=
			\sum \limits_{k= j+1}^\infty \text{Pr} (\mathfrak{N}^\nu(t) = k \big|
			\mathfrak{N}^\nu(0) = 1 ) \\ \notag
			& = 1- \sum \limits_{k= 1}^j \text{Pr} (\mathfrak{N}^\nu(t) = k \big|
			\mathfrak{N}^\nu(0) = 1 )\\
			& = 1- \sum \limits_{k= 1}^j  \sum \limits_{l= 1}^k \binom{k-1}{l-1}
			(-1)^{l-1} E_{\nu,1} (-\lambda l t^\nu). \notag
		\end{align}
		This implies that the $j$th waiting time $\mathfrak{W}_j^\nu$  has distribution
		\[ %begin{equation}
			f_{\mathfrak{W}_j^\nu} (t)= \sum \limits_{k= 1}^j
			\sum \limits_{l= 1}^k \binom{k-1}{l-1} (-1)^{l-1} (\lambda l) t^{\nu-1}
			E_{\nu, \nu} (-\lambda l t^\nu), \quad t>0, \: \nu \in (0,1].
		\] %end{equation}
		Integrating the preceding equation, we get
		\begin{align}  \notag
			\int \limits_0^\infty f_{\mathfrak{\mathfrak{W}}_j^\nu} (t) dt &=
			\sum \limits_{k= 1}^j   \sum \limits_{l= 1}^k \binom{k-1}{l-1} (-1)^{l-1} \\ \notag
			\notag
			&= \sum \limits_{k= 1}^j  \sum \limits_{l= 0}^{k-1} \binom{k-1}{l} (-1)^l= \sum
			\limits_{k= 1}^j  (1-1)^{k-1}=1. \notag
		\end{align}
		The non-negativity of $f_{\mathfrak{W}_j^\nu} (t)$ follows from the non-negativity of
		$ \mathfrak{p}_k^\nu(t)$ (see \ocite{pol}),  and the last line of \eqref{eqprob}
		is a monotone increasing function of  $t$. To see this, we can write
		\begin{align}  \notag
			&1- \sum \limits_{k= 1}^j  \sum \limits_{l= 1}^k \binom{k-1}{l-1}
			(-1)^{l-1} E_{\nu,1} (-\lambda l t^\nu) \\  \notag
			&= 1-\sum \limits_{k= 1}^j  \frak{p}_k^\nu(t) \notag \\ \notag
			& =  1- \sum \limits_{k= 1}^j  \sum \limits_{l= 1}^k \binom{k-1}{l-1} (-1)^{l-1}
			\text{Pr} (T_l >t) \\ \notag
			&= \sum \limits_{k= 1}^j  \sum \limits_{l= 1}^k \binom{k-1}{l-1} (-1)^{l-1}
			\text{Pr} (T_l < t ). \notag
		\end{align}
		Indeed, $f_{\mathfrak{W}_j^\nu} (t) $ is a probability density function.
		Note also that  $f_{\mathfrak{W}_j^\nu} (t)$ has the following integral
		representation:
		\[ %begin{equation}
			f_{\mathfrak{W}_j^\nu} (t) = \frac{1}{t} \int_0^\infty e^{-\xi}
			\sum \limits_{k= 1}^j   \sum \limits_{l= 1}^k \binom{k-1}{l-1} (-1)^{l-1}
			g( l \lambda t / \xi) d\xi,
		\] %end{equation}
		where $g(\eta) =\sin (\nu \pi) / [ \pi( \eta^\nu+ \eta^{-\nu} + 2\cos (\nu \pi ) ) ]$
		(see \ocite{repin}).
		We now show that if the sojourn times are distributed as in \eqref{interbirth},
		the cumulative  distribution function $\text{Pr} ( \mathfrak{W}_j^\nu \leq t) $
		of the waiting or birth time equals the right-hand side of  \eqref{eq2a}. When $j=1$, we get
		\[ %begin{equation}
			\text{Pr} (\mathfrak{W}_1^\nu \leq t)= \text{Pr} ( T_1^\nu \leq t) =
			1 - E_{\nu, 1} (- \lambda t^\nu )= 1- \mathfrak{p}_1^\nu (t).
		\] %end{equation}
		In the succeeding calculations, we use the following identities (see page 26 of \ocite{podlubny}):
		\[ %begin{equation}
			\int_0^t E_{\nu, 1} (- j \lambda (t-u)^\nu)  u^{\nu-1} E_{\nu, \nu} (-\lambda l u^\nu) du
			= \frac{jE_{\nu, \nu +1} (- j \lambda t^\nu)-l E_{\nu,\nu +1}(-l\lambda t^\nu)}{j-l}t^{\nu }
		\]%end{equation}
		and
		\[ %begin{equation}
			E_{\nu, \nu + 1} (\xi) = \frac{ E_{\nu, 1} (\xi) -1}{\xi}, \qquad l \leq j.
		\] %end{equation}

		%\noindent When $j=2$,
		\noindent Now,
		\begin{align} \notag
			\text{Pr} (\mathfrak{W}_2^\nu  \leq t) &= \int_0^t \text{Pr} \lbrace
			T_1^\nu + T_2^\nu  \leq t \big| T_1^\nu = u   \rbrace  dF_{T_1^\nu} (u)\\ \notag
			&= \int_0^t [ 1- E_{\nu, 1} (- 2 \lambda (t-u)^\nu) ]  \lambda u^{\nu-1}
			E_{\nu, \nu} (- \lambda u^\nu ) du \\\notag
			&= 1 - E_{\nu, 1} (- \lambda t^\nu) - [  2 \lambda t^\nu  E_{\nu, \nu +1}
			(- 2 \lambda t^\nu) - t^\nu E_{\nu, \nu +1} (-  \lambda t^\nu) ] \\\notag
			&= 1 - E_{\nu, 1} (- \lambda t^\nu) - [ E_{\nu, 1} (- \lambda t^\nu) -
			E_{\nu, 1} (-2\lambda t^\nu) ] \\ \notag
			&= 1 - 2E_{\nu, 1} (- \lambda t^\nu) + E_{\nu, 1} (-2\lambda t^\nu)  \\ \notag
			&= 1-\sum_{k= 1}^2 \frak{p}_k^\nu (t), \notag
		\end{align}
		and in general, we can show by induction that
		\begin{align} \notag
			\text{Pr}& (\mathfrak{W}_j^\nu   \leq t) \\ \notag
			= {} & \int_0^t \text{Pr} \lbrace   \mathfrak{W}_j^\nu  \leq t \big|
			\mathfrak{W}_{j-1}^\nu  = u   \rbrace  dF_{\mathfrak{W}_{j-1}^\nu} (u) \notag \\ \notag
			= {} & \int_0^t [ 1- E_{\nu, 1} (- j \lambda (t-u)^\nu) ]
			f_{\mathfrak{W}_{j-1}^\nu} (u) du \\ \notag
			= {} & \int_0^t [ 1- E_{\nu, 1} (- j \lambda (t-u)^\nu) ]
			\sum \limits_{k= 1}^{j-1}  \sum \limits_{l= 1}^k \binom{k-1}{l-1} (-1)^{l-1}
			(\lambda l) u^{\nu-1} E_{\nu, \nu} (-\lambda l u^\nu) du \\ \notag
			= {} & \sum \limits_{k= 1}^{j-1}  \sum \limits_{l= 1}^k \binom{k-1}{l-1} (-1)^{l-1}
			[1-E_{\nu,1} (-\lambda l t^\nu)] \\\notag
			& -   \sum \limits_{k= 1}^{j-1}  \sum \limits_{l= 1}^k \binom{k-1}{l-1} (-1)^{l-1}
			\lambda l \int_0^t E_{\nu, 1} (- j \lambda (t-u)^\nu)  u^{\nu-1} E_{\nu, \nu}
			(-\lambda l u^\nu) du\\ \notag
			= {} & \sum \limits_{k= 1}^{j-1}  \sum \limits_{l= 1}^k \binom{k-1}{l-1} (-1)^{l-1}
			[1-E_{\nu,1} (-\lambda l t^\nu)] \\ \notag
			&- \sum \limits_{k= 1}^{j-1}  \sum \limits_{l= 1}^k \binom{k-1}{l-1} (-1)^{l-1}
			\frac{l}{j-l}[E_{\nu,1} (-\lambda l t^\nu) -E_{\nu,1} (-\lambda j t^\nu) ]  \\ \notag
			= {} & 1- \sum \limits_{k= 1}^{j-1}  \sum \limits_{l= 1}^k \binom{k-1}{l-1} (-1)^{l-1}
			\left(\frac{j}{j-l}E_{\nu,1} (-\lambda l t^\nu) - \frac{l}{j-l}E_{\nu,1} (-\lambda jt^\nu)
			\right)  \\ \notag
		        = {} & 1- \sum \limits_{k= 1}^{j-1}  \sum \limits_{l= 1}^k \binom{k-1}{l-1} (-1)^{l-1}
			\frac{j}{j-l}E_{\nu,1} (-\lambda l t^\nu) \\ \notag
		        & +  \sum \limits_{k= 1}^{j-1}  \sum \limits_{l= 1}^k \binom{k-1}{l-1} (-1)^{l-1}
			\frac{l}{j-l}E_{\nu,1} (-\lambda jt^\nu).  \\ \notag
		\end{align}
		\noindent Using the formulas on page 3 of \ocite{gar80}, we have
		\begin{align}  \notag
			\sum \limits_{k= 1}^{j-1}  \sum \limits_{l= 1}^k \binom{k-1}{l-1} (-1)^{l-1} \frac{l}{j-l}
			&= \sum \limits_{l= 1}^{j-1} (-1)^{l-1} \frac{l}{j-l}  \sum \limits_{k =l}^{j-1}
			\binom{k-1}{l-1}   \\ \notag
			&= \sum \limits_{l= 1}^{j-1} (-1)^{l-1} \frac{l}{j-l}  \sum \limits_{k =0}^{j-1-l}
			\binom{k + l-1}{l-1}   \\ \notag
			&= \sum \limits_{l= 1}^{j-1} (-1)^{l-1} \frac{l}{j-l} \binom{j-1}{l}   \\ \notag
			&= \sum \limits_{l= 1}^{j-1} (-1)^{l-1} \frac{l}{j-l}\frac{(j-1)!}{l!(j-l-1)!}  \\ \notag
			&= \sum \limits_{l= 1}^{j-1} \frac{(j-1)!}{(l-1)!(j-l)!}  \\ \notag
			&= \sum \limits_{l= 0}^{j-2}  (-1)^l \binom{j-1}{l} =
			(-1)^{j-2}, \\ \notag
		\end{align}
		because
		\begin{align}
			\sum_{l=0}^{j-2} (-1)^l \binom{j-1}{l} = \sum_{l=0}^{j-1} (-1)^l \binom{j-1}{l} -
			(-1)^{j-1} \binom{j-1}{j-1}. \notag
		\end{align}
		\noindent Hence,
		\begin{align}
			\text{Pr}& (\mathfrak{W}_j^\nu \leq t) \\       \notag
			= {} & 1- \sum \limits_{k= 1}^{j-1}  \sum \limits_{l= 1}^k \binom{k-1}{l-1} (-1)^{l-1}
			\frac{j}{j-l}E_{\nu,1} (-\lambda l t^\nu) - (-1)^{j-1} E_{\nu,1} (-\lambda j t^\nu)
			\\ \notag
			= {} & 1- \sum \limits_{k= 1}^{j-1}  \sum \limits_{l= 1}^k \binom{k-1}{l-1} (-1)^{l-1}
			\frac{j}{j-l}E_{\nu,1} (-\lambda l t^\nu) \\
			& - \left( \sum \limits_{l= 1}^j \binom{j-1}{l-1}
			(-1)^{l-1}  E_{\nu,1} (-\lambda l t^\nu) \notag
			-  \sum \limits_{l= 1}^{j-1} \binom{j-1}{l-1}
			(-1)^{l-1}  E_{\nu,1} (-\lambda l t^\nu)  \right) \\ \notag
			= {} & 1- \sum \limits_{k= 1}^{j-1}  \sum \limits_{l= 1}^k \binom{k-1}{l-1}
			(-1)^{l-1}	\frac{j}{j-l}E_{\nu,1} (-\lambda l t^\nu) \\
		        & - \left( \sum \limits_{l= 1}^j \binom{j-1}{l-1}
			(-1)^{l-1}  E_{\nu,1} (-\lambda l t^\nu) \notag
		        -  \sum \limits_{l= 1}^{j-1} \binom{j-1}{l}
			(-1)^{l-1}\frac{l}{j-l}  E_{\nu,1} (-\lambda l t^\nu)  \right) \\ \notag
 %		\end{align}
 %\begin{align}
			= {} & 1- \sum \limits_{k= 1}^{j-1}  \sum \limits_{l= 1}^k \binom{k-1}{l-1}
			(-1)^{l-1} \frac{j}{j-l}E_{\nu,1} (-\lambda l t^\nu) \\
		        & - \left( \sum \limits_{l= 1}^j \binom{j-1}{l-1}
			(-1)^{l-1}  E_{\nu,1} (-\lambda l t^\nu) \notag
			-  \sum \limits_{l= 1}^{j-1} (-1)^{l-1}\frac{l}{j-l}  E_{\nu,1} (-\lambda l t^\nu)
			\sum \limits_{k= l}^{j-1} \binom{k-1}{l-1}  \right) \\ \notag
			= {} & 1- \sum \limits_{k= 1}^{j-1}  \sum \limits_{l= 1}^k \binom{k-1}{l-1}
			(-1)^{l-1} \frac{j}{j-l}E_{\nu,1} (-\lambda l t^\nu) \\
			& - \left( \sum \limits_{l= 1}^j \binom{j-1}{l-1}
			(-1)^{l-1}  E_{\nu,1} (-\lambda l t^\nu) \notag
			-  \sum \limits_{l= k}^{j-1} \sum \limits_{l= 1}^k  \binom{k-1}{l-1}
			(-1)^{l-1} \frac{l}{j-l} E_{\nu,1} (-\lambda l t^\nu)    \right) \\ \notag
			= {} & 1- \sum \limits_{k= 1}^{j-1}  \sum \limits_{l= 1}^k \binom{k-1}{l-1}
			(-1)^{l-1}E_{\nu,1} (-\lambda l t^\nu) - \sum \limits_{l= 1}^j \binom{j-1}{l-1}
			(-1)^{l-1}  E_{\nu,1} (-\lambda l t^\nu) \\ \notag
			= {} & 1-\sum \limits_{k= 1}^j \frak{p}_k^\nu (t), \quad 1\leq  k < j, \notag
		\end{align}
as the second summation (in the preceding equal sign) simply corresponds to $k=j$. Hence, equality \eqref{eq2a} is attained.  Again,  the transition from the third equality to the fourth equality above uses formula (0.15.1) on page 3 of \ocite{gar80}, i.e.,
\[
\sum_{k=l}^{j-1} \binom{k-1}{l-1} = \sum_{k=0}^{j-l-1}
\binom{k+l-1}{l-1} = \binom{j-1}{l}.
\]
Notice that when $\nu=1$, we get
	$\text{Pr}(\mathfrak{W}_j \leq t)=(1-e^{-\lambda t})^j$  which corresponds to the birth time distribution of the classical Yule
		process. Moreover, equality \eqref{eq2b} can be straightforwardly evaluated as
		\begin{align} \notag
			\text{Pr} ( \mathfrak{W}_{j-1}^\nu \leq t) - \text{Pr} ( \mathfrak{W}_j^\nu \leq t) =
			{} & \left(1- \sum \limits_{k= 1}^{j-1}  \sum \limits_{l= 1}^k \binom{k-1}{l-1}
			(-1)^{l-1} E_{\nu,1} (-\lambda l t^\nu) \right) \\ \notag
			& - \left( 1- \sum \limits_{k= 1}^j  \sum \limits_{l= 1}^k \binom{k-1}{l-1}
			(-1)^{l-1} E_{\nu,1} (-\lambda l t^\nu) \right) \\ \notag
			= {} & \sum \limits_{l= 1}^j \binom{j-1}{l-1} (-1)^{l-1} E_{\nu,1}
			(-\lambda l t^\nu) \\ \notag
			= {} & \frak{p}_j^\nu (t). \notag
		\end{align}
		In addition,  the Laplace transform of the probability density $f_{T_i^\nu} (t)$ is
		\[ %begin{equation}
			\int_0^\infty e^{-zt} f_{T_i^\nu} (t) dt = \frac{i \lambda}{i \lambda + z^\nu}.
		\] %end{equation}
		This suggests that the distribution (eqn \eqref{interbirth}) leads to the following known
		mixture or structural representation (see \ocite{cuw10}) of the inter-birth times as
		\[ %begin{equation}
			T_i^\nu \stackrel{d}{=}  V_i^{1/\nu} S_{\nu},
		\] %end{equation}
		where $V_i$ has the exponential distribution with parameter $ i \lambda$, i.e.,
		\begin{equation} \label{expd}
			f_{V_i}(v)=i \lambda e^{-i \lambda v}, \quad  v >0,
		\end{equation}
		and is independent of the positive $L\acute{e}vy$ or $\nu$-stable distributed random variable
		$S_{\nu}$  having the Laplace transform of the  density function $ e^{-z^\nu}$.
		This also suggests that the $\kappa$-$th$ fractional moment of the $i$th inter-birth time is given by
		\[ %begin{equation}
			\mathbf{E} \left[ T_i^\nu \right]^\kappa =\frac{\pi \Gamma (1 + \kappa)}{(i \lambda)^\kappa \Gamma
			(\kappa / \nu) \sin (\pi \kappa/ \nu ) \Gamma (1 - \kappa)}, \quad 0< \kappa <\nu,
		\] %end{equation}
		which further implies that the $\kappa$-$th$ fractional moment of the $j$th wait or birth time is
		\[ %begin{equation}
			\mathbf{E} \left[ \mathfrak{W}_j^\nu \right]^\kappa = \frac{\pi \Gamma (1 + \kappa)}{\lambda^\kappa
			\Gamma (\kappa / \nu) \sin (\pi \kappa/ \nu ) \Gamma (1 - \kappa)} \sum \limits_{k= 1}^j
			\sum \limits_{l= 1}^k \binom{k-1}{l-1} (-1)^{l-1} \left( \frac{1}{l^\kappa} \right),
		\] %end{equation}
		where $0<\kappa<\nu$.

%\newpage
%	\section{Parameter Estimation}

	\section{Sample paths of fYp}

		From Sections 3 and 4, it is now straightforward to simulate a trajectory of a fYp.
       % using the mixture representation of the random inter-birth times.
        However, we only propose the two simplest algorithms on how to generate a sample path of the fYp as the others
		follow.  In particular,  the random-rate representation (Representation A, Theorem \ref{repra})
		yields the algorithm below.

\vspace{0.05in}

		\noindent ALGORITHM 1:
\vspace{-0.05in}
		\begin{description}
			\item i) Generate $\Xi$ from the Wright distribution $W_{-\nu, 1-\nu} (-\xi)$, and obtain $\xi$.
			\item ii) Simulate a classical Yule process with birth rate $\lambda \xi$.
			\item iii) Stretch the time scale to $t^\nu$.
		\end{description}

		A simpler way to generate a realization of fYp with $n$ births is to directly exploit the known
		birth and/or sojourn time distributions as follows: Generate
		$V_i$ from the exponential distribution in \eqref{expd} with parameter $ i \lambda$, and
		$S_\nu$ from the strictly positive stable distribution with parameter $\nu$.

\vspace{0.05in}

%\pagebreak
		\noindent ALGORITHM 2:
\vspace{-0.05in}
		\begin{description}
			\item i) Let $i=1$ and $\mathfrak{N}^\nu (0)=1$.
			\item ii) Simulate $T_i^\nu = V_i^{1/ \nu}S_\nu$, and let $\mathfrak{W}_i^\nu =
				T_1^\nu + T_2^\nu + \cdots +  T_i^\nu$ .
			\item iii) $\mathfrak{N}^\nu (\mathfrak{W}_i^\nu)=i+1$, and $i=i+1$.
			\item iv) Repeat ii-iii  for $i=2,\ldots,n-1$.
		\end{description}

We now use the algorithms above to highlight some unique properties of the fractional Yule process that are related to its true mean given in (\ref{mun}). Figure \ref{Figure1} below shows  both Yp and fYp as jump processes of size 1  in the time interval $(0, 5)$  with $\nu=0.5$, and $\lambda=1$. Using the same  set of parameters, Figure   \ref{Figure2}  displays sample trajectories of a different/independent fYp and Yp which model a binary-split growth process.  An important attribute that can be directly observed from these two graphs is that on the average, fYp  grows more rapidly than the classical Yp  shortly after it starts.

%\bigskip

%\bigskip

		\begin{figure}[h!t!b!p!]
			\centering
			\includegraphics[height=2.5in, width=4.7in]{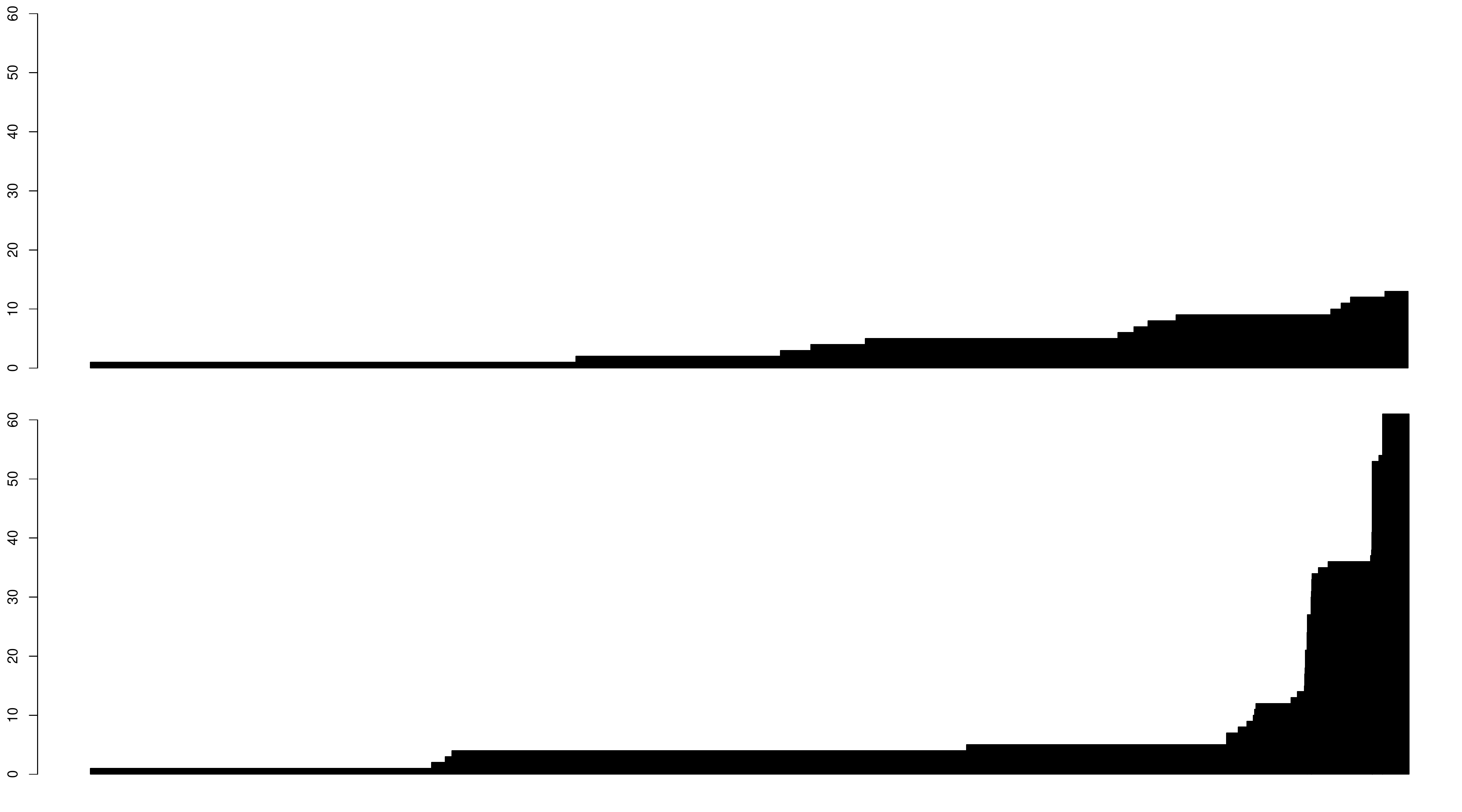}
			\caption{\emph{Sample trajectories of  the standard Yule process  (top)  and  the fractional Yule process (bottom) in the interval $(0, 5)$ with parameters $ (\nu, \lambda)=(0.5, 1)$. }}
			\renewcommand\belowcaptionskip{0pt}
			\label{Figure1}
		\end{figure}

		\begin{figure}[h!t!b!p!]
			\centering
			\includegraphics[height=2.5in, width=5in]{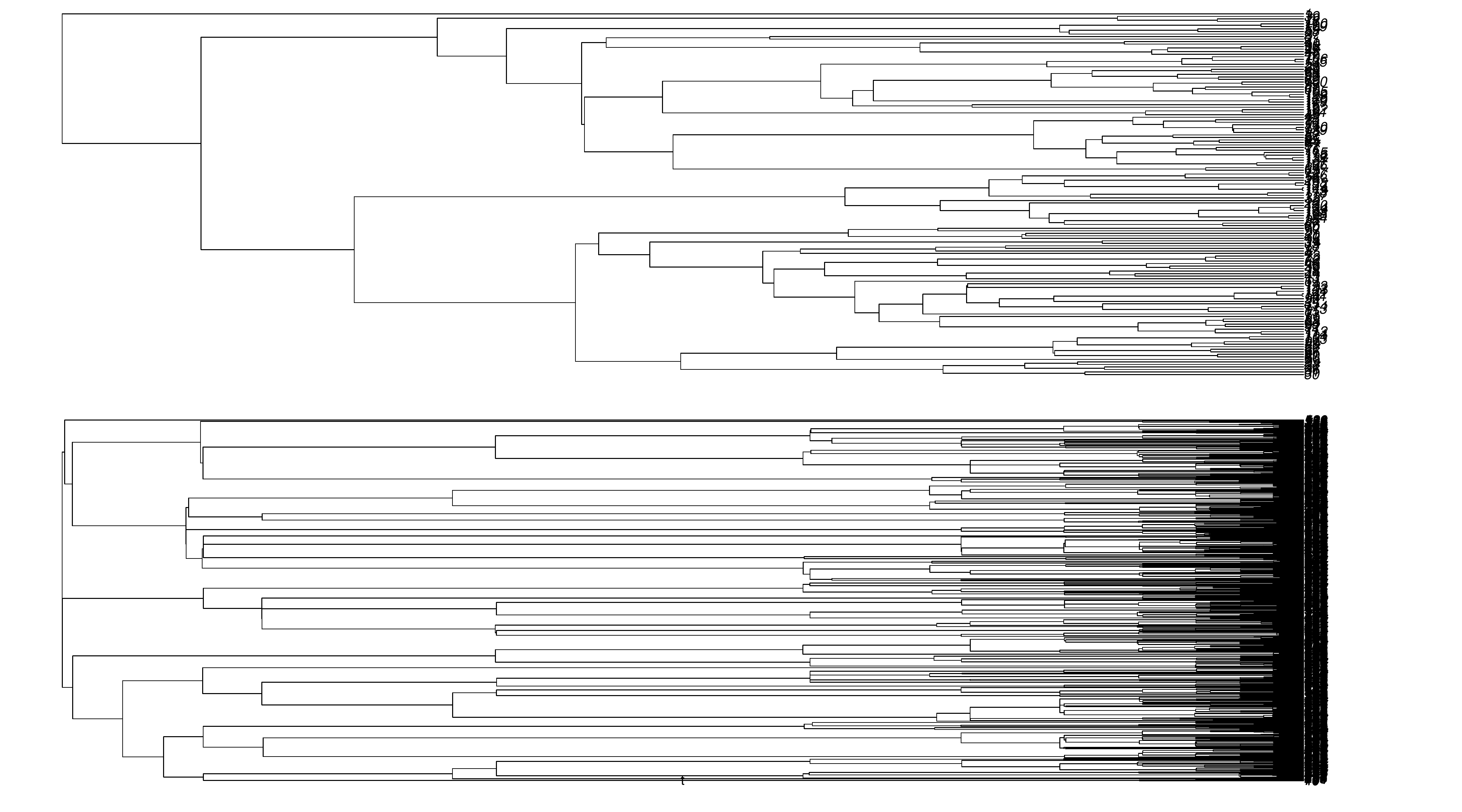}
			\caption{\emph{Sample paths of the pure linear birth process  (top) and  the fractional Yule
				process (bottom) in the interval $(0, 5)$ with parameters   $ (\nu, \lambda)=(0.5, 1).$ }}
			 \renewcommand\belowcaptionskip{0pt}
			\label{Figure2}
		\end{figure}

\begin{figure}[h!t!b!p!]
			\centering
			\includegraphics[height=3in, width=5in]{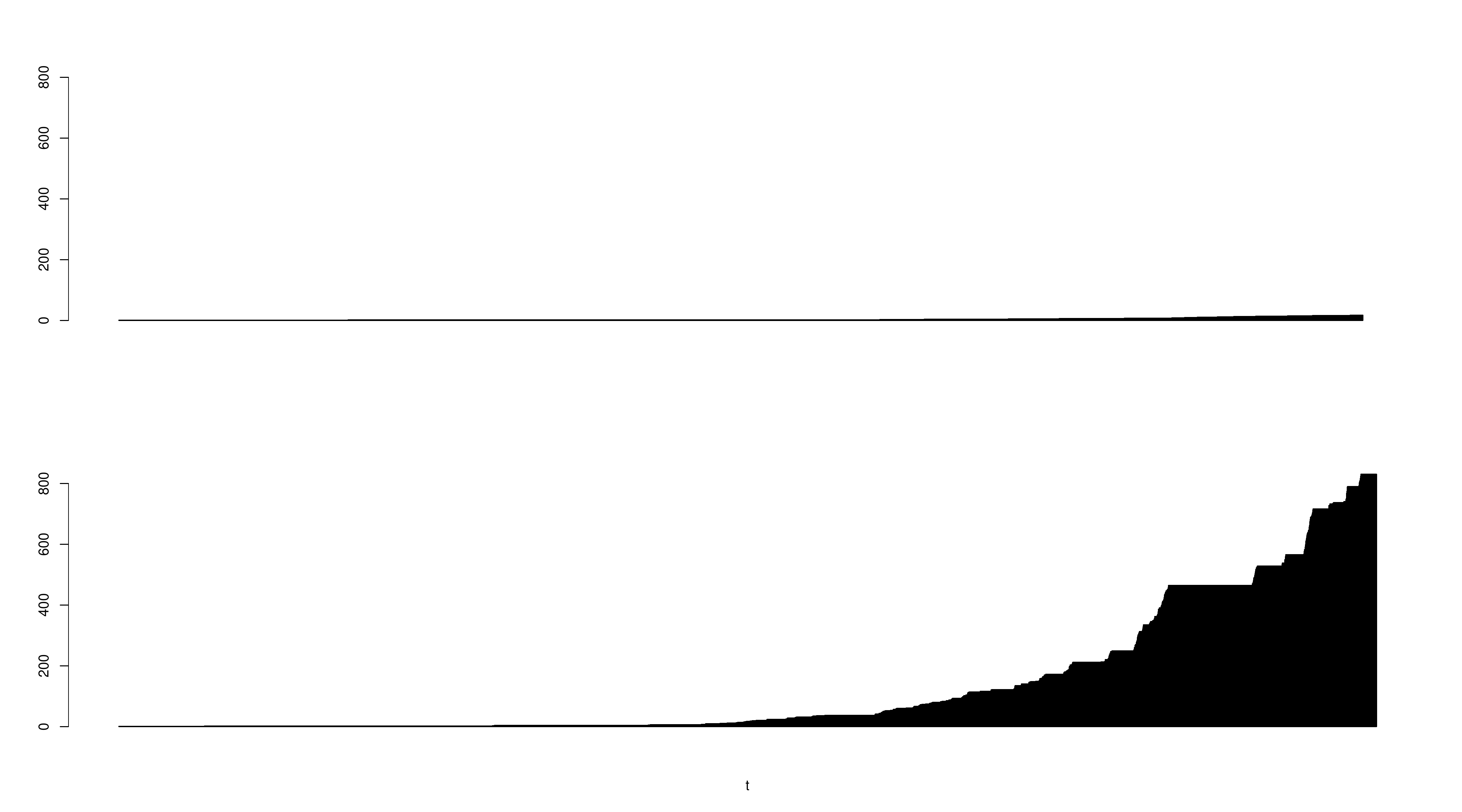}
			\caption{\emph{Sample paths of the classical Yule process  (top) and  the fractional Yule
				process (bottom) in the interval $(0, 5)$ with parameters   $(\nu, \lambda)=(0.25, 1)$. }}
			 \renewcommand\belowcaptionskip{0pt}
			\label{Figure3}
		\end{figure}

\pagebreak
In addition, a more specific  characteristic of  fYp is illustrated  in  Figure  \ref{Figure3}. The particular realization of fYp below used the parameter values  $\nu=0.25$, $\lambda=1$, and is observed in the time interval $(0, 5)$ . It clearly suggests that fYp is more explosive than Yp when $\nu \to 0$. In general, the plots strongly validate the  plausibility of fYp to model exploding and strictly growing processes. Note also that Representation A implies that the interaction between  the random rate  and time stretching of the classical Yule process can rapidly  speed up or slow down fYp at any given time instance.

\begin{comment}

		\begin{figure}[h!t!b!p!]
			\centering
			\includegraphics[height=3.1in, width=5in]{path1c.pdf}
			\caption{\emph{Sample paths of classical Yule process  (top) and  fractional Yule
				process (bottom)  in the interval $(0, 5)$ with parameters   $ (\nu, \lambda)=(0.99, 1).$ }}
			 \renewcommand\belowcaptionskip{0pt}
			\label{Figure4}
		\end{figure}
\end{comment}

\pagebreak

	\section{Method-of-Moments (MoM) estimation}

		We now propose a method-of-moments  estimation procedure for the parameters $\nu$ and $\lambda$ to
		make fYp usable in practice. In this procedure, we assume that a particular realization
		or complete history of the process is observed until the population is $n$, i.e., there are $n$
		births. We then
		attempt to use all the available data from the observed sample path of the fractional Yule process.

		In particular, we use all the available inter-birth or sojourn times of the observed sample
		trajectory of the fractional Yule process. A  direct way of estimating the parameters is to use the fractional moment estimators as follows: Choose constants $\kappa_m < \nu,  m=1,2$, and solve for the estimates $\hat{\lambda}$ and $\hat{\nu}$ using the equations
		\[ %begin{equation}
			\frac{\sum_{i=1}^n \left[ T_i^\nu \right]^{\kappa_m} }{n}= \frac{\pi
			\Gamma (1 + \kappa_m)}{ \hat{\lambda}^{\kappa_m} \Gamma (\kappa_m / \hat{\nu} ) \sin (\pi \kappa_m /
			\hat{\nu} ) \Gamma (1 - \kappa_m)}\frac{\sum_{i=1}^n 1/i^{\kappa_m}}{n} , \quad m=1,2.
		\] %end{equation}
		Another approach is to use the first two  integer-order moments  of the log-transformed sojourn
		times (see \ocite{cuw10}) which are
		\[ %begin{equation}
			\mathbf{E} \ln \left[ T_i^\nu \right] = \frac{- \ln (i \lambda)}{\nu} - \gamma,
		\] %end{equation}
		and
		\[ %begin{equation}
			\mathbf{E} \ln \left[ T_i^\nu \right]^2 = \pi^2 \left( \frac{1}{3 \nu^2} -
			\frac{1}{6} \right) +  \left( \frac{ \ln (i \lambda)}{\nu} + \gamma \right)^2.
		\] %end{equation}
		This further suggests that the parameter estimates can be computed using the two equations:
		\[ %begin{equation}
			\frac{\sum_{i=1}^n \ln \left[ T_i^\nu \right] }{n}= \frac{- \sum_{i=1}^n \ln
			(i \lambda)}{\nu  n} - \gamma,
		\] %end{equation}
		and
		\[ %begin{equation}
			\frac{\sum_{i=1}^n  \left( \ln \left[ T_i^\nu \right]\right)^2 }{n}=
			\pi^2 \left( \frac{1}{3 \nu^2} - \frac{1}{6} \right) + \frac{1}{n}\sum_{i=1}^n
			\left( \frac{ \ln (i \lambda)}{\nu} + \gamma \right)^2,
		\] %end{equation}
		where $\gamma \cong 0.577215664901532 $ is the Euler--Mascheroni  constant. A major advantage of this procedure over other moment estimators is that it does not require selection of  constants a priori to calculate the parameter estimates.  Note  also that the maximum likelihood estimators  are more challenging to compute due to the required evaluation of the 		 Mittag--Leffler function.

		In addition, we  tested our parameter estimation procedure. In doing so, we generated 10 random samples of inter-birth times
		of size 10000 each for  $\nu=0.1 + 0.1m$, 	$m=0, \ldots, 9$ and $\lambda=0.2,10$.  For each simulated data set, we computed the estimates
		using the first $n$ observations in the set with $n=100, 1000$, and $10000$.
		The tables below show the simulation results for a single run, which further indicate that the proposed procedure
		performs relatively well as the sample sizes increase.  Please note that in many applications (e.g., internet traffic), the typical number of observations is at least of the order of millions.  These estimates could also serve as good starting values of an iterative estimation procedure.

		\begin{table}[h!t!b!p!]
			\caption{\emph{Parameter estimates $(\widehat{\nu}, \widehat{\lambda})$ for fYp  with
				$\nu=0.1(0.1)1$ and $\lambda=0.2$. }} \centerline {
			\begin{tabular*}{4.5in}{@{\extracolsep{\fill}}|c||ccc|} %|l@{\hspace{0.01in}}||c@{\hspace{0.01in}}c@{\hspace{0.01in}}|}%{lp{0.01\textwidth}||cp{0.1\textwidth}cp{0.01\textwidth}|}
				\hline
				 & $n=100$ & $n=1000$ & $n=10000$\\
				\hline \hline
				$(\nu=0.1, \lambda=0.2)$ & (0.095, 0.198)  & (0.096, 0.185)  & (0.100, 0.205)   \\
				$(\nu=0.2, \lambda=0.2)$ & (0.228, 0.249)  & (0.193, 0.189)  & (0.199, 0.193)  \\
				$(\nu=0.3, \lambda=0.2)$ & (0.283, 0.185) &  (0.292, 0.193)  & (0.303, 0.228) \\
				$(\nu=0.4, \lambda=0.2)$ & (0.381, 0.178) &  (0.407, 0.218)  & (0.402, 0.209)  \\
				$(\nu=0.5, \lambda=0.2)$ & (0.481, 0.212)  & (0.501, 0.197)  & (0.500, 0.197)   \\
				$(\nu=0.6, \lambda=0.2)$ & (0.599, 0.211)  & (0.602, 0.186)  & (0.595, 0.186)  \\
				$(\nu=0.7, \lambda=0.2)$ & (0.759, 0.257) &  (0.728, 0.250)  & (0.700, 0.198) \\
				$(\nu=0.8, \lambda=0.2)$ & (0.818, 0.220) &  (0.819, 0.229)  & (0.803, 0.204)   \\
				$(\nu=0.9, \lambda=0.2)$ & (0.850, 0.193)  & (0.899, 0.211)  & (0.907, 0.215)   \\
				$(\nu=1.0, \lambda=0.2)$ & (0.977, 0.183)  & (0.991, 0.199)  & (0.999, 0.202)  \\
				\hline
			\end{tabular*}
			}
		\end{table}

		\begin{table}[h!t!b!p!]
			\caption{\emph{Parameter estimates $(\widehat{\nu}, \widehat{\lambda})$ for fYp  with
			$\nu=0.1(0.1)1$ and $\lambda=10$. }} \centerline {
			 \begin{tabular*}{4.5in}{@{\extracolsep{\fill}}|c||ccc|}%|l@{\hspace{0.01in}}||c@{\hspace{0.01in}}c@{\hspace{0.01in}}|}%{lp{0.01\textwidth}||cp{0.1\textwidth}cp{0.01\textwidth}|}
				\hline
				 & $n=100$ & $n=1000$ & $n=10000$ \\
				\hline \hline
				$(\nu=0.1, \lambda=10)$ & (0.107, 13.067)  & (0.101, 10.599)  & (0.101, 10.730) \\
				$(\nu=0.2, \lambda=10)$ & (0.203, 10.737)  & (0.206, 12.384)  & (0.201, 10.555)  \\
				$(\nu=0.3, \lambda=10)$ & (0.299, 11.027) & (0.297, 9.359)  & (0.295, 8.593) \\
				$(\nu=0.4, \lambda=10)$ & (0.391, 7.598)  & (0.396, 8.899)  & (0.397, 9.086)   \\
				$(\nu=0.5, \lambda=10)$ & (0.517, 10.939) & (0.509, 11.428) & (0.501, 10.269)   \\
				$(\nu=0.6, \lambda=10)$ & (0.630, 11.379) & (0.586, 8.308)  & (0.597, 9.162)  \\
				$(\nu=0.7, \lambda=10)$ & (0.716, 12.413) & (0.699, 10.634)  & (0.710, 11.679) \\
				$(\nu=0.8, \lambda=10)$ & (0.782, 8.713)  & (0.786, 8.186)  & (0.804, 10.498)   \\
				$(\nu=0.9, \lambda=10)$ & (0.919, 11.429) & (0.899, 9.043)  & (0.897, 9.684)   \\
				$(\nu=1.0, \lambda=10)$ & (0.969, 8.712)  & (1.000, 10.427)  & (1.001, 10.434)  \\
				\hline
			\end{tabular*}
			}
		\end{table}

	%\pagebreak

	\section{Concluding remarks}

		We have derived one-dimensional representations of the fractional Yule process, which led to algorithms  for simulating its sample paths. These representations are also
		necessary in understanding  the properties of fYp further.  We have derived the birth and inter-birth or sojourn
		time distributions, which are of Mittag--Leffler type. The structural representation of the random
		sojourn time also led to an algorithm for simulating sample trajectories of the fYp.  We have proposed
		an  estimation procedure using the  moments of the log-transformed inter-birth times, which
		performed satisfactorily especially for larger sample sizes.

		Although some properties of fYp have already been studied, there are  still a lot of open problems
		that need to be figured out. For instance, understanding fYp in more depth and the construction of more efficient estimators like the maximum
		likelihood would be worth pursuing in the future. Also, the  application of fYp in practice  particularly in biology and/or network traffic is still in progress.

%\begin{acknowledgements}
%\end{acknowledgements}

\bibliographystyle{spphys}       % APS-like style for physics

\end{document}